\documentclass{article}
\usepackage{fancyhdr}
\usepackage{amscd}
\usepackage{graphicx}
\usepackage{amsmath}
\usepackage{amssymb}
\usepackage{amsthm}
\usepackage{mathrsfs}
\newtheorem{dfn}{Definition}[section]
\newtheorem{thm}[dfn]{Theorem}

\newtheorem{rem}[dfn]{Remark}

\numberwithin{equation}{section}

\begin{document}

\title{A refined asymptotic behavior of traveling wave solutions for degenerate nonlinear parabolic equations}

\author{Yu Ichida\thanks{Graduate School of Science and Technology, Meiji University, 1-1-1, Higashimita Tama-ku Kawasaki Kanagawa 214-8571, Japan, {\tt ce196001@meiji.ac.jp}}, Kaname Matsue\thanks{Institute of Mathematics for Industry, Kyushu University, Fukuoka 819-0395, Japan {\tt kmatsue@imi.kyushu-u.ac.jp}} $^{,}$ \footnote{International Institute for Carbon-Neutral Energy Research, Kyushu University, Fukuoka 819-0395, Japan} $^{,}$ \footnote{Center for Research and Development Strategy, Japan Science and Technology Agency (JST-CRDS), Tokyo 102-0076, Japan}, Takashi Okuda Sakamoto\thanks{Graduate School of Science and Technology, Meiji University, 1-1-1, Higashimita Tama-ku Kawasaki Kanagawa 214-8571, Japan, {\tt sakamoto@meiji.ac.jp}}
}
%\date{}
\maketitle

\begin{abstract}
 In this paper,  we consider the asymptotic behavior of traveling wave solutions of the degenerate nonlinear parabolic equation: $u_{t}=u^{p}(u_{xx}+u)-\delta u\, (\delta = 0 \, \mbox{or} \,  1)$ for $\xi \equiv x - ct \to - \infty$ with $c>0$.
We give a refined one of them, which was not obtain in the preceding work \cite{DNPE}, by an appropriate asymptotic study and properties of the Lambert $W$ function.
\end{abstract} 


\section{Introduction}
In this paper, we consider the degenerate nonlinear parabolic equation
\begin{equation}
u_{t}=u^{p}(u_{xx}+u)-\delta u,\quad t>0,\quad x\in \mathbb{R},
\label{eq:dnla1}
\end{equation}
where $\delta=0$ or $1$, $p\in 2\mathbb{N}$.

When $\delta=0$, this equation arises in the modeling of heat combustion, solar flares in astrophysics, plane curve evolution problems and the resistive diffusion of a force-free magnetic field in a plasma confined between two walls (see \cite{Anada, Low1, Low2, Poon} and references therein).
Also, there are many studies on blow-up solution to \eqref{eq:dnla1} (for instance, see \cite{Anada, Poon} and references therein).

On the other hand, the equation \eqref{eq:dnla1} with $\delta=1$ can be obtained by transforming solution of \eqref{eq:dnla1} with $\delta=0$ (see \cite{Poon}).
In \cite{Poon},  the traveling wave solutions of \eqref{eq:dnla1} play important roles. 
More precisely, the lower bound of the blow-up rate is obtained by means of the traveling wave solutions of \eqref{eq:dnla1} under either the Dirichlet boundary condition or the periodic boundary condition in the case that $\delta=1$ and $x$ is restricted to $x\in (-L,L)$. 

In addition, the traveling wave solutions are not only upper (or lower) solutions as discussed in \cite{Poon} but also the entire solutions of the equation. 
These facts motivate us to study detailed information of the traveling wave solutions to \eqref{eq:dnla1}.

In order to consider the traveling waves of \eqref{eq:dnla1}, we introduce the following change of variables:
\[ u(t, x)=\phi(\xi), \quad \xi=x-ct, \quad c>0. \]
The equation of $\phi(\xi)$ solving \eqref{eq:dnla1} is then reduced to
\begin{equation}
-c\phi'=\phi^{p}\phi''+\phi^{p+1}-\delta\phi, \quad \xi \in \mathbb{R}, \quad '=\dfrac{d}{d\xi}, \nonumber
\end{equation}
equivalently 
\begin{equation}
\left \{
\begin{array}{l}
\phi' = \psi, \\
\psi' = -c\phi^{-p}\psi-\phi+\delta\phi^{-p+1},
\end{array}
\right. \label{eq:dnla3} 
\end{equation}
where $\delta = 0$ or $1$. 

In \cite{DNPE}, a result on the whole dynamics on the phase space $\mathbb{R}^{2}$ including infinity
generated by the two-dimensional ordinary differential equation (ODE for short) \eqref{eq:dnla3} is obtained by applying the dynamical system approach and the Poincar\'e compactification (for instance, see  \cite{FAL} for the details of the Poincar\'e compactification).
Further, connecting orbits on it are focused and several results on the existence of (weak) traveling wave solutions are given.
The following theorem is one of main results obtained in \cite{DNPE}.

\begin{thm}[\cite{DNPE}, Theorem 3]
 \label{th:dnla1}
Assume that $p\in 2\mathbb{N}$ and $\delta=1$.
Then, for a given positive constant $c$, the equation \eqref{eq:dnla1} has a family of traveling wave solutions (which correspond to a family of the orbits of \eqref{eq:dnla3}).
Each traveling wave solution $u(t,x) = \phi(\xi)$ satisfies the following.
\begin{itemize}
\item [$\bullet$]
$\left\{ \begin{array}{ll}
%\displaystyle \lim_{\xi \to -\infty} u(\xi) = 0, & 
%\quad \displaystyle \lim_{\xi \to +\infty} u(\xi) = 1,\\
%\displaystyle\lim_{\xi \to -\infty} u'(\xi) = 0,& 
%\quad \displaystyle\lim_{\xi \to +\infty} u'(\xi) =  0.
\displaystyle \lim_{\xi \to -\infty} \phi(\xi) = 0, & 
\quad \displaystyle \lim_{\xi \to +\infty} \phi(\xi) = 1,\\
\displaystyle\lim_{\xi \to -\infty} \phi'(\xi) = 0,& 
\quad \displaystyle\lim_{\xi \to +\infty} \phi'(\xi) =  0.
\end{array} \right.$

\item [$\bullet$]
$\phi(\xi)>0$ holds for $\xi\in \mathbb{R}$.

\end{itemize}
\end{thm}

Figure \ref{fig:dnla1} shows dynamics on the Poincar\'e disk of \eqref{eq:dnla3} (see \cite{FAL} for the definition of the Poincar\'e disk).  In addition, the asymptotic behavior of the traveling wave solutions (obtained in Theorem \ref{th:dnla1}) for $\xi \to +\infty$ is also given in \cite{DNPE}, while the asymptotic behavior as $\xi \to -\infty$ is not obtained there.

%\begin{figure}[t]
%\centering
%\includegraphics[scale=0.12]{dnla-keijyou.pdf}
%\caption{Schematic picture of the traveling wave solution on $\xi \in \mathbb{R}$ obtained in Theorem\ref{th:dnla1} in the case that $D=c^{2}-4p>0$.}
%\label{fig:dnla-kei}
%\end{figure}

In this paper, we give a refined asymptotic behavior of the traveling wave solutions, which contributes to extraction of their characteristic nature.
The main theorem of this paper is the following.

\begin{thm}
\label{th:dnla2}
%Assume that $p\in 2\mathbb{N}$ and $\delta=1$.
%Then, for a given positive constant $c$, the equation \eqref{eq:dnla1} has a family of traveling wave solutions (which correspond to a family of the orbits of \eqref{eq:dnla3}).
%Each traveling wave solution $u(\xi)$ satisfies the following.
%\begin{itemize}
%\item [$\bullet$]
%$\left\{ \begin{array}{ll}
%\displaystyle \lim_{\xi \to -\infty} u(\xi) = 0, & 
%\quad \displaystyle \lim_{\xi \to +\infty} u(\xi) = 1,\\
%\displaystyle\lim_{\xi \to -\infty} u'(\xi) = 0,& 
%\quad \displaystyle\lim_{\xi \to +\infty} u'(\xi) =  0.
%\end{array} \right.$

%\item [$\bullet$]
%$u(\xi)>0$ holds for $\xi\in \mathbb{R}$.

%\end{itemize}
The asymptotic behavior of $\phi(\xi)$ obtained in Theorem \ref{th:dnla1} as $\xi \to -\infty$ is 
\begin{equation*}
\phi(\xi) \sim \left(  \dfrac{\mu c^{2}}{ \mu (c^{2}+1)-e^{-\frac{p}{c}\xi}} \right)^{\frac{1}{p}}, \quad {\rm as} \quad \xi \to -\infty,
 %\label{eq:dnla-asy}
\end{equation*}
where $\mu<0$ is a constant that depends on the initial state $\phi_{0} = \phi(0)$.
\end{thm}

During our proof of the theorem, we see that {\em the Lambert $W$ function}
plays a key role in describing the asymptotic behavior. 
Evaluation of integrals including the Lambert $W$ function is necessary to obtain the asymptotic behavior in the present form.
Our argument here is based on an asymptotic study of solutions in the different form from that provided in e.g. \cite{DNPE, Matsue2}, which can be applied to asymptotic analysis towards further applications in various phenomena including their numerical calculations.

%In order to prove Theorem\ref{th:dnla2}, it is necessary to prove the asymptotic behavior of $u(\xi)$ for $\xi \to -\infty$.
%It should be noted again that such conditions are satisfied and the existence of traveling wave solutions was proved in \cite{DNPE}.

\section{Preliminaries}
In this section, we partially reproduce calculations in \cite{DNPE} for the readers' convenience.

First, we study the dynamics near bounded equilibria of \eqref{eq:dnla3}.
If $\delta=1$ and $p$ is even, then \eqref{eq:dnla3} has the equilibria $\pm E_{\delta} : (\phi,\psi)=(\pm 1,0)$.
Let $J_{1}$ be the Jacobian matrix of the vector field \eqref{eq:dnla3} at $E_{\delta}$.
Then, the behavior of the solution around $E_{\delta}$ is different by the sign of $D=c^{2}-4p$.
For instance, the matrix $J_{1}$ has the real distinct eigenvalues if $D>0$ and other cases can be concluded similarly.
In addition, if $c>0$, then the real part of all eigenvalues of $J_{1}$ are negative.
Therefore, we determine that the the equilibria $\pm E_{\delta} : (\phi,\psi)=(\pm1,0)$ are sink.

Second, in order to study the dynamics of \eqref{eq:dnla3} on the Poincar\'e disk, we desingularize it by the 
time-scale desingularization 
\begin{equation} 
ds/d\xi = \{\phi(\xi)\}^{-p}\quad  {\rm for} \quad  p \in 2\mathbb{N}. 
\label{eq:dnla4}
\end{equation}
Since $p$ is assumed to be even, the direction of the time direction does not change via this desingularization.
Then we have
\begin{equation}
\left \{
\begin{array}{l}
\phi' = \phi^{p}\psi, \\
\psi' = -c\psi-\phi^{p+1}+\delta\phi,
\end{array}
\right. \quad \left(~'{\mbox{}}~ = \dfrac{d \mbox{}}{ds}\right), 
\label{eq:dnla5}
\end{equation}
where $\delta=0$ or $\delta=1$.

It should be noted that the time-scale desingularization \eqref{eq:dnla4} is simply the multiplication of $\phi^p$ to the vector field. Then, except the singularity $\{\phi=0\}$, the solution curves of the system (vector field) remain the same but are parameterized differently (see also Section 7.7 of \cite{CK}). 
%Still, we refer to Section 7.7 of \cite{CK} and references therein for the analytical treatments of desingularization with the time rescaling. 
%In what follows, we use the similar time rescaling (re-parameterization of the solution curves) repeatedly to desingularize the vector fields.

The system \eqref{eq:dnla5} has the equilibrium $E_{O}:(\phi,\psi)=(0,0)$.
When $\delta=1$, the Jacobian matrix of the vector field \eqref{eq:dnla5} at $E_{O}$ is 
\[ E_{O} : \left(\begin{array}{cc}
0 & 0 \\
1 & -c
\end{array}
\right). \]
It has the real distinct eigenvalues $0$ and $-c$. 
The eigenvectors corresponding to each eigenvalue are 
\[ \mathbf{v}_{1}=\left(\begin{array}{cc}
c \\
1
\end{array}
\right),\quad \mathbf{v}_{2}=\left(\begin{array}{cc}
0 \\
1
\end{array}
\right). \]
We set a matrix $T$ as $T=(\mathbf{v}_{1},\mathbf{v}_{2})$.
Then we obtain
\begin{align*}
\left(\begin{array}{cc}
\phi' \\
\psi'
\end{array}
\right) &= \left(\begin{array}{cc}
0 & 0 \\
1 & -c
\end{array}
\right)\left(\begin{array}{cc}
\phi \\
\psi
\end{array}
\right)+\left(\begin{array}{cc}
\phi^{p}\psi \\
-\phi^{p+1}
\end{array}
\right) \\
&= T\left(\begin{array}{cc}
0 & 0 \\
0 & -c
\end{array}
\right)T^{-1}\left(\begin{array}{cc}
\phi \\
\psi
\end{array}
\right)+\left(\begin{array}{cc}
\phi^{p}\psi \\
-\phi^{p+1}
\end{array}
\right).
\end{align*}
%By multiplying $T^{-1}$ to both sides from left, 
%We 
Let 
$\left(\begin{array}{cc}
\tilde{\phi} \\
\tilde{\psi}
\end{array}
\right)=T^{-1}\left(\begin{array}{cc}
\phi \\
\psi
\end{array}
\right)$.
%Then, we can see the following equations:
%\[ \phi^{p}\psi=c^{p}\tilde{\phi}^{p+1}+c^{p}\tilde{\phi}^{p}\tilde{\psi}, \quad -\phi^{p+1}=-c^{p+1}\tilde{\phi}^{p+1}. \]
%Therefore, we can obtain the following equation:
We then obtain the following system:
\begin{equation}
\left \{
\begin{array}{l}
\tilde{\phi}'=c^{p-1}\tilde{\phi}^{p+1}+c^{p-1}\tilde{\phi}^{p}\tilde{\psi}, \\
\tilde{\psi}'=-c\tilde{\psi}-c^{p-1}\tilde{\phi}^{p+1}-c^{p-1}\tilde{\phi}^{p}\tilde{\psi}-c^{p+1}\tilde{\phi}^{p+1}.
\end{array}
\right.  
\label{eq:dnla6}
\end{equation}
The center manifold theory (e.g. \cite{carr}) is applicable to study the dynamics of \eqref{eq:dnla6}.
It implies that there exists a function $h(\phi)$ satisfying 
\[ h(0)=\dfrac{dh}{d\tilde{\phi}}(0)=0 \]
such that the center manifold of $E_O$ for \eqref{eq:dnla6} is locally represented as $\{(\tilde{\phi},\tilde{\psi}) \,|\, \tilde{\psi}(s)=h(\tilde{\phi}(s))\}$. 
Differentiating it with respect to $s$, we have
\begin{align*}
-ch(\tilde{\phi})-c^{p-1}\tilde{\phi}^{p+1}-c^{p-1}\tilde{\phi}^{p}h(\tilde{\phi})-c^{p+1}\tilde{\phi}^{p+1}
\\=\dfrac{dh}{d\tilde{\phi}} \left( c^{p-1}\tilde{\phi}^{p+1}+c^{p-1}\tilde{\phi}^{p}h(\tilde{\phi}) \right). 
\end{align*}
Then we obtain the approximation of the (graph of) center manifold as follows:
\begin{equation}
\left\{ (\tilde{\phi}, \tilde{\psi}) \,|\, \tilde{\psi}=-c^{p-2}(c^{2}+1)\tilde{\phi}(s)^{p+1}+O(\tilde{\phi}^{p+2}) \right\}. 
\label{eq:dnla7} 
\end{equation}
Therefore, the dynamics of \eqref{eq:dnla6} near $E_O$ is topologically equivalent to the dynamics of the following equation:
\begin{equation*}
\tilde{\phi}'(s)=c^{p-1}\tilde{\phi}^{p+1}-c^{2p-3}(c^{2}+1)\tilde{\phi}^{2p+1}.
%\label{eq:cm4} 
\end{equation*}
%Note that $\tilde{\phi}$ and $\tilde{\psi}$ are $\tilde{\phi}=\phi/c$, $\tilde{\psi}=\psi-\phi/c$.
We conclude that the approximation of the (graph of) center manifold are 
\begin{equation}
\phi'(s) =\phi^{p+1}/c- [(c^{2}+1)\phi^{2p+1}]/c^{3}
\label{eq:dnla8} 
\end{equation}
and the dynamics of \eqref{eq:dnla5} near $E_O$ is topologically equivalent to the dynamics of the following equation:
\begin{equation}
\psi(s) =\phi/c-[(c^{2}+1)\phi^{p+1}]/c^{3}. \nonumber
%\label{eq:dnla9} 
\end{equation}

\begin{figure}[t]
\centering
\includegraphics[scale=0.23]{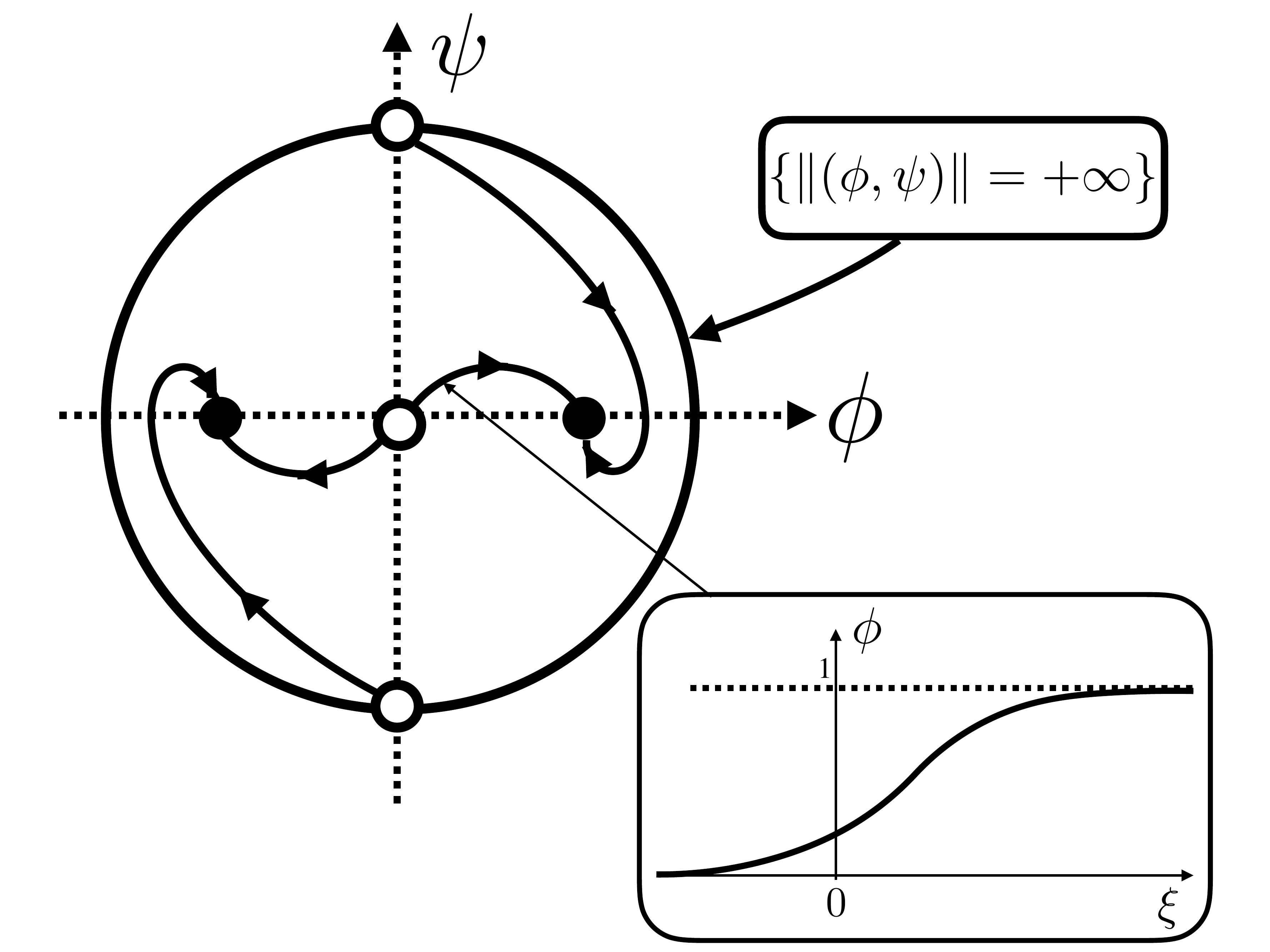}
\caption{Schematic picture of the dynamics on the Poincar\'e disk and corresponding traveling wave solutions in the case that $\delta$ is $\delta=1$ and $p$ is even with $D=c^{2}-4p<0$ and $c>0$.}
\label{fig:dnla1}
\end{figure}

Finally, we obtain the dynamics on the Poincar\'e disk in the case that $p$ is even (see Figure\ref{fig:dnla1}).
This argument indicates that the asymptotic behavior of $\phi$ through the present system is calculated as a function of $s$ and that an additional asymptotic study is required to obtain the behavior of $\phi$ in terms of the original frame coordinate $\xi$.

\begin{rem}[\cite{DNPE}, Remark 1]
\label{rem:dnla1}
In Figure \ref{fig:dnla1}, we need to be careful about the handling of the point $E_{O}$.
When we consider the parameter $s$ on the disk, $E_{O}$ is the equilibrium of \eqref{eq:dnla5}.
However, $E_{O}$ is a point on the line $\{\phi=0\}$ with singularity about the parameter $\xi$.
We see that $d \phi/d \psi$ takes the same values on the vector fields defined by \eqref{eq:dnla5} and \eqref{eq:dnla3} except the singularity $\{\phi=0\}$.
If the trajectories start the equilibrium $E_{O}$ about the parameter $s$, then they start from the point $E_{O}$ about $\xi$. 
\end{rem}

\section{Proof of Theorem \ref{th:dnla2}}
The proof is divided into four steps. 
In Step I, we derive an ODE describing the behavior of $s$ with respect to $\xi$. 
It turns out to contain {\em the Lambert $W$ function}.
In Step II, we confirm that $\xi(s) \to -\infty$ as $s \to -\infty$, which is used for the direct derivation of $\phi(\xi)$ in the asymptotic sense. 
Step III is devoted to obtain the relationship between $\phi$ and $\xi$. 
According to preceding studies such as \cite{DNPE, Matsue2}, the asymptotic behavior of $\phi(\xi)$ can be obtained {\em in the composite form} $\phi(s(\xi))$, which can require multiple integrations of differential equations.
Except special cases, lengthy calculations are necessary towards an explicit and meaningful expression of the targeting asymptotics.
Instead, we directly derive the relationship of $\phi$ to $\xi$ without solving the ODE obtained in Step (I) and calculate the asymptotic behavior of the function $\xi(\phi)$ as $\phi \to 0$ associated with the center manifold \eqref{eq:dnla7}, which works well even if integrands include the Lambert W function.
We finally obtain the asymptotic behavior of $\phi(\xi)$ in Step IV via inverse function arguments.
\\
\begin{rem}
\label{rem:dnla2}
{\em The Lambert $W$ function} $y = W(x)$ is defined as the inverse function of $x=ye^{y}$. 
We easily see the following properties which we shall use below:
\begin{itemize}
\item $W(x) >0$ for $x>0$;
\item $W(x) < \log x$ for $x>e$. 
\end{itemize}
See e.g. \cite{Lam} and references therein for further properties.
\end{rem}
\mbox{}\\
{\textbf{Proof.}}
{\bf (I)}:
First we set
\[ w(s):=\phi(s)^{-p}>0. \]
With the aid of (\ref{eq:dnla8}), we have
\begin{equation}
w'(s) =-p/c + [ p(c^{2}+1) / c^{3} w] =A+Bw^{-1},  
\label{eq:dnla10}
\end{equation}
where 
%\[ A=-p/c\neq 0 \quad \mbox{ and } \quad B=[p(c^{2}+1)/c^{3}]\neq 0.\]
\[  A=-p/c < 0 \quad \mbox{ and } \quad B=[p(c^{2}+1)/c^{3}] > 0.\]
The solution of \eqref{eq:dnla10} satisfies the following.
\[ \left| 1+Aw/B\right| e^{-\left(\frac{A}{B}w+1\right)}= \left|A/B\right| e^{-\frac{A^{2}}{B}s-\frac{A^{2}C_{1}+B}{B}}  \]
with a constant $C_{1}$.
Since the dynamics of $\phi(s)$ near $0$ (i.e., $\phi(s) \approx 0$) is of our interest, 
we may assume that $w(s)$ is sufficiently large, which implies that $[A/B]w+1<0$.
Then we have
\[ -\left(1+Aw/B\right) e^{-\left(\frac{A}{B}w+1\right)}= -Ae^{-\frac{A^{2}}{B}s-\frac{A^{2}C_{1}+B}{B}} /B. \]
By using $w=\phi(s)^{-p}$ and the Lambert $W$ function, we obtain
\[
\phi(s)= \left[-B\left\{ W\left( E(s) \right)+1\right\}/A \right]^{-\frac{1}{p}}, 
\]
where $E(s)=-[A/B]e^{-\frac{A^{2}}{B}s-\frac{A^{2}C_{1}+B}{B}}$.
We consequently have
\begin{equation}
\dfrac{ds}{d\xi} =\phi^{-p}
=-B\left\{ W\left(E(s)\right)+1\right\} /A.
\label{eq:dnla11}
\end{equation}
%This yields
%\[
%\xi+C_{2}=\int \left[ -B \left\{ W\left(E(s)\right)+1\right\}/A \right]^{-1}ds
%\]
%with a constant $C_{2}$.
\\
{\bf (II)}: We shall prove
\[ \xi(s) \to -\infty \quad \mbox{as} \quad s \to - \infty.\]
We note that $E(\cdot)$ is positive on $\mathbb{R}$ and hence $W(E(s))>0$ holds for $s \in \mathbb{R}$.
Integrating \eqref{eq:dnla11} on $(-\infty,0]$, we have
\[ \xi(0) - \xi_{-} = \int_{-\infty}^{0} \left[ -B \left\{ W\left(E(s)\right)+1\right\}/A \right]^{-1}ds, \]
where
\[ \xi_- = \lim_{s \to -\infty} \xi(s).\]
Without loss of generality, we may set $\xi(0) = 0$. 

By using properties of the Lambert $W$ function, for a negative constant $s_*$ satisfying $|s_* | \gg 1 $,
we have
\begin{align*}
-\xi_{-} 
&=   \int_{-\infty}^{0} \left[ -B \left\{ W\left(E(s)\right)+1\right\}/A \right]^{-1}ds\\
&>   - \frac{A}{B} \Bigg [\int_{s_*}^{0}   \left\{ W\left(E(s)\right)+1\right\}^{-1}ds  +  \int_{-\infty}^{s_*} \{\log(E(s))+1\}^{-1} ds \Bigg ]\\
& > - \frac{A}{B}   \int_{-\infty}^{s_*} \left[ \log \left( -[A/B]e^{-\frac{A^{2}}{B}s-\frac{A^{2}C_{1}+B}{B}}\right )+1\right ]^{-1} ds\\
& = - \frac{A}{B}   \int_{-\infty}^{s_*} [ \, \log( -[A/B])  - ({A^{2}}s+{A^{2}C_{1}+B} )/B+1  \, ] ^{-1} ds\\
&=  - \lim_{s \to -\infty}\frac{A}{B} \bigg[ \, (- B/A^2) \log| (-A^2/B)s + C_2 |  \, \bigg]_{s}^{s_\ast},
\end{align*}
where
\[ C_2 = \log(-A/B) - (A^2C_1+B)/B + 1. \]
Since $A<0<B$ holds, we have
\begin{align*}
 -\xi_{-} & > (1/A) \log\{ (-A^2/B)s_*+C_2\} \\
& \quad  + (-1/A) \lim_{s \to -\infty} \log \{(-A^2/B)s+ C_2 \}  = + \infty.
\end{align*}
Therefore the asymptotic behavior of $\phi(s)$ as $s\to -\infty$ is equivalent to that of $\phi(\xi)$ as $\xi \to -\infty$.
\\
\mbox{}\\
{\bf (III)}: Next, we represent $\xi$ as a function of $\phi$.
We rewrite \eqref{eq:dnla11} as
\[ \dfrac{d\xi}{ds}=\phi^{p}.\]
Using \eqref{eq:dnla8}, we obtain
\begin{align*}
\xi+C_{3} &=\int  \{\phi(s)\}^{p} ds = \int \phi^{p}\dfrac{ds}{d\phi}d\phi 
\\
%&= \int \phi^{p}\dfrac{ds}{d\phi}d\phi 
%\\
&= \int \phi^{p} \left(\dfrac{1}{c}\phi^{p+1}-\dfrac{c^{2}+1}{c^{3}}\phi^{2p+1}\right)^{-1} d\phi 
\\
&=\int \dfrac{c^{3}}{\phi\{c^{2}-(c^{2}+1)\phi^{p}\}}d\phi
\end{align*}
with a constant $C_{3}$.
Introducing $\varphi=\phi^{p}$, we further have
\begin{align*}
\xi+C_{3} &=\int \dfrac{c^{3}}{\phi\{c^{2}-(c^{2}+1)\phi^{p}\}}d\phi 
\\
&=\dfrac{c^{3}}{p}  \int \dfrac{1}{ \varphi \{c^{2}-(c^{2}+1)\varphi\}}  d\varphi
\\
&= \dfrac{c^{3}}{p}  \int \left\{\dfrac{1}{c^{2}}\dfrac{1}{\varphi}+\left(1+\dfrac{1}{c^{2}}\right)\dfrac{1}{ c^{2}-(c^{2}+1)\varphi }\right\} d\varphi
\\
%&= \dfrac{c}{p} \int \dfrac{1}{\varphi}d\varphi -\dfrac{c^{3}}{p}\left(1+\dfrac{1}{c^{2}}\right) \int \dfrac{1}{ (c^{2}+1)\varphi-c^{2}} d\varphi
%\\
&=\dfrac{c}{p}\log \left| \dfrac{\varphi}{(c^{2}+1)\varphi-c^{2}} \right|
\\
&=\dfrac{c}{p} \log \left| \dfrac{\phi^{p}}{(c^{2}+1)\phi^{p}-c^{2}} \right|.
\end{align*}
Then the constant $C_{3}$ is given by 
\[ C_{3}=\dfrac{c}{p} \log \left| \dfrac{\phi_{0}^{p}}{(c^{2}+1)\phi_{0}^{p}-c^{2}} \right|, \]
where $\phi(0)=\phi_{0}$.
Moreover, it holds that $C_{3}<0$ regardless of the value of $c$, provided $\phi_0 \ll 1$.
Indeed, it holds that 
\[ 0< \phi_{0} <\left( c^{2}/ [ c^{2}+2 ] \right)^{\frac{1}{p}} \qquad (0 < \phi_0 \ll 1). \]
%Here we note again that we consider the dynamics near $\phi \approx 0$.
\\
%Moreover, it holds that if $C_3 < 0$ then $c>0$ and vice versa. } \\
{\bf (IV)}:  Finally, we aim to represent $\phi$ as a function of $\xi$.
As mentioned above, we obtain 
\[ \xi+\dfrac{c}{p} \log \left| \dfrac{\phi_{0}^{p}}{(c^{2}+1)\phi_{0}^{p}-c^{2}} \right|= \dfrac{c}{p} \log \left| \dfrac{\phi^{p}}{(c^{2}+1)\phi^{p}-c^{2}} \right|. \]
This yields
%\[ \dfrac{\phi^{p}}{(c^{2}+1)\phi^{p}-c^{2}} = \pm \left| \dfrac{\phi_{0}^{p}}{(c^{2}+1)\phi_{0}^{p}-c^{2}} \right| e^{\frac{p}{c}\xi} =\mu e^{\frac{p}{c}\xi}
% \]
%with a constant $\mu$.
\[ \dfrac{\phi^{p}}{(c^{2}+1)\phi^{p}-c^{2}} = \pm \left| \dfrac{\phi_{0}^{p}}{(c^{2}+1)\phi_{0}^{p}-c^{2}} \right| e^{\frac{p}{c}\xi} .
 \]
%Therefore, we have
%\[
%\phi^{p} = \dfrac{\mu e^{\frac{p}{c}\xi}c^{2}}{ \mu e^{\frac{p}{c}\xi}(c^{2}+1)-1}.
%\]
Therefore, we have
\[
\phi^{p} = \dfrac{\mu e^{\frac{p}{c}\xi}c^{2}}{ \mu e^{\frac{p}{c}\xi}(c^{2}+1)-1}, \quad 
\mu=\pm \left| \dfrac{\phi_{0}^{p}}{(c^{2}+1)\phi_{0}^{p}-c^{2}} \right|.
\]

If $\mu>0$, there exists a finite  $\xi $ such that $\mu e^{\frac{p}{c}\xi}(c^{2}+1)-1=0$ holds.
%In other words, it occurs singularity.
%However, this is contradicted by the fact that there is no minimum value for $\xi$.
However, as in Theorem \ref{th:dnla1}, the traveling wave solutions $\phi(\xi)$ that correspond to the connecting orbits between
$E_O$ and $E_{\delta}$ have no singularities for $\xi \in \mathbb{R}$.
Therefore, $\mu$ must be negative. 
This yields 
\[ \dfrac{\mu e^{\frac{p}{c}\xi}c^{2}}{ \mu e^{\frac{p}{c}\xi}(c^{2}+1)-1}>0
\quad \mbox{with} \quad \mu= - \left| \dfrac{\phi_{0}^{p}}{(c^{2}+1)\phi_{0}^{p}-c^{2}} \right|. \]
Since $p$ is even, we obtain the following.
\[ 
\phi(\xi) =\left(  \dfrac{\mu c^{2}}{ \mu (c^{2}+1)-e^{-\frac{p}{c}\xi}} \right)^{\frac{1}{p}} \to 0, \quad {\rm as} \quad \xi \to -\infty,
 \]
 where $\mu<0$ is the constant that depends on the initial state $\phi_{0}$.
\qed

\section{Conclusion}
In this paper, we give a refined asymptotic behavior of the traveling wave solutions of \eqref{eq:dnla1} as $\xi \to -\infty$.
As shown in Step III of the proof, the present result is obtained by considering the asymptotic behavior of $\xi(\phi)$ without taking the relationship between $\xi$ and $s$ into account.
This is a key idea to get over the difficulties of treatment of the Lambert $W$ function to obtain the asymptotic behavior for $u(t,x) = \phi(\xi)$.
We expect that our approach can be applied to the asymptotic behavior of typical solutions as well as that of singular solutions.

\section*{Acknowledgments}
KM was partially supported by World Premier International Research Center Initiative (WPI), Ministry of Education, Culture, Sports, Science and Technology (MEXT), Japan and the grant-in-aid for young scientists No.~17K14235, Japan Society for the Promotion of Science.

\end{document}